\documentclass{amsart}
\newtheorem{theorem}{Th\'eor\`eme}[section]
\newtheorem{lemma}[theorem]{Lemme}
\newtheorem{prop}[theorem]{Proposition}
\newtheorem{corollary}[theorem]{Corollaire}
\theoremstyle{definition}

\theoremstyle{remark}
\newtheorem{remark}[theorem]{Remarque}

\numberwithin{equation}{section}

\def\H{\mathbb H}

\def\Z{\mathbb Z}

\def\SS{\mathcal S}

\def\NN{\mathcal N}
\def\HH{\mathcal H}
\def\LL{\mathcal L}

\def\n{\noindent}

\newcommand{\pt}{ \,^{p}\!\tau  }

\newcommand{\ilm}{  j_{!*}{\mathcal L}}
\newcommand{\ph}[2]{ \,^p\!{\mathcal H}^{#1}(#2)}

\newcommand{\im}{ \hbox{\rm Im} }

\title{\large   De la Puret\'e locale \`a la d\'ecomposition }
\author{Fouad El Zein}
\address{Institut de Math\'ematiques de Jussieu,
Paris, France}
 \email{elzein@math.jussieu.fr}
 \author{ D\~ung Tr\'ang L\^e}
\address{Universit\'e d'Aix-Marseille 
LATP, UMR-CNRS 7353
39, rue Joliot-Curie
F-13453 Marseille Cedex 13}
 \email{ledt@ictp.it}

\date{Janvier  2013}
\begin{document}
\maketitle 

\begin{abstract} Le th\'eor\`eme de d\'ecomposition se d\'eduit de la puret\'e locale.
\end{abstract}

\vskip.2in
\centerline {\large From Local purity  to Decomposition}

\begin{abstract}
The decomposition theorem is deduced from  local purity.
\end{abstract}
\section{Introduction}
 Cette note fait suite \`a la note \cite{EL} o\`u nous avons
expliqu\'e comment le th\'eor\`eme de d\'ecomposition en dehors
d'un point
implique  un th\'eor\`eme  de puret\'e locale en ce point dans le cadre analytique, similaire au th\'eor\`eme
de puret\'e locale de Deligne-Gabber de \cite{DG} dans le cadre alg\'ebrique. Ici nous montrons comment
ce th\'eor\`eme de puret\'e locale en un point sert \`a \'etendre le th\'eor\`eme de d\'ecomposition
en ce point, ce qui permet d'\'etablir simultan\'ement par r\'ecurrence le th\'eor\`eme  de puret\'e locale 
et le th\'eor\`eme de d\'ecomposition.

 Nous reprenons les notations de \cite{EL}. 
 Soit  $f: X \to V$ un  morphisme  projectif de 
 vari\'et\'es alg\'ebriques complexes, $\tilde \LL$ 
une variation de structures de Hodge  (VSH) polaris\'ee  sur un ouvert $ \Omega$ lisse de 
$X$, $j: \Omega \to X$, $\LL := \tilde \LL [m]$ le complexe 
de cohomologie r\'eduite \`a $\tilde \LL$ en degr\'e  $ - m $ o\`u $ m  $
est la dimension de $X$, et $j_{!*} \LL$ l'extension  interm\'ediaire de  $\LL$. Dans 
la note pr\'ec\'edente \cite{EL} nous avons d\'eduit la puret\'e locale en un point $v$ de $V$, 
de la d\'ecomposition de  $Rf_*j_{!*}\LL$ sur $V-v$ . Il s'agit \`a pr\'esent d'\'etendre la 
d\'ecomposition  en $v$. 

Cette notion consiste en r\'ealit\'e en deux d\'ecompositions de nature bien distinctes. L'une consiste \`a 
d\'ecomposer une cohomologie perverse 
 $ \ph{i}{Rf_*(j_{!*}\LL)}$  en une somme directe d'extensions 
interm\'ediaires canoniques 
qui utilise la th\'eorie de Hodge.

L'autre consiste \`a d\'ecomposer le complexe en une somme directe
de  complexes \'el\'ementaires se r\'eduisant \`a ses propres cohomologies  perverses  d\'ecal\'ees 
en degr\'e, ce  qui d\'ecoule de  la d\'eg\'en\'erescense de la suite spectrale de Leray 
perverse et qui se d\'eduit  d'un r\'esultat de  type Lefschetz \cite{SS}, \`a savoir: le  
 cup-produit  it\'er\'e avec la classe $\eta$ d'une section 
   hyperplane, induit des isomorphismes
 de (complexes de) faisceaux de  cohomologie perverse
$
\eta^i: {^p\!\HH}^{-i}(Rf_*j_{!*}\LL)
\xrightarrow{\sim}{^p\!\HH}^i(Rf_*j_{!*}\LL).
$ 

La d\'emontration  proc\`ede par une r\'ecurrence sur les strates d'une stratification convenable  du morphisme $f$. Le pas de r\'ecurrence consiste \`a appliquer le th\'eor\`eme  de puret\'e locale \`a une section normale \`a une strate  de $V$ en un point  de la strate. 
Pour poursuivre le raisonnement  de r\'ecurrence et utiliser  chaque fois une SHM sur le compl\'ementaire d'un diviseur \`a croisements normaux (DCN), on doit pr\'eparer la vari\'et\'e $X$ \`a l'aide 
d'une d\'esingularisation adapt\'ee \`a $\LL$ et \`a une stratification de Thom-Whitney de $f$, 
ce qui ne change pas la port\'ee du r\'esultat puisque l'on peut toujours se r\'eduire \`a ce cas. 

En fait, la d\'ecomposition refl\`ete  les propri\'et\'es topologiques des morphismes projectifs en  g\'eom\'etrie alg\'ebrique complexe et il est agr\'eable de lier ces r\'esultats \`a des fibrations topologiques induites sur les strates d'une  stratification de Thom-Whitney sur la base. De plus, ayant choisi d'utiliser 
des complexes logarithmiques, on s'int\'eresse particuli\`erement \`a des fibrations par des DCN sur les strates au sens de la note pr\'ec\'edente  (voir D\'efinition 3.1 \cite{EL}).

On {\it d\'ebute la r\'ecurrence} sur l'ouvert $U$
form\'e par la r\'eunion des grandes strates lisses  de $V$ tel que  la restriction de $f$ au-dessus de $U$ soit  lisse et 
propre. On suppose par construction que:   $f^{-1}(U)$  est le 
compl\'ementaire d'un DCN, le compl\'ementaire $Y = X -  \Omega$ est un DCN et que les singularit\'es de $\LL$  forment au-dessus de $U$
 un DCN horizontal $Y_U := Y \cap f^{-1}(U) $ relatif sur $U$.
 
Alors la famille de  cohomologie d'intersection des fibres forme une VSH polaris\'ee. Dans ce cas,  les   faisceaux image-directes  sup\'erieures   $R^if_*j_{!*}\LL$
 classiques  coincident  avec les cohomologies perverses: $\ph{i}{Rf_*
(j_{!*}\LL)} = {\HH}^{i}(Rf_*(j_{!*}\LL))$, le  th\'eor\`eme de Lefschetz difficile  
s'applique, et par cons\'equent
 les r\'esultats de \cite {SS} aussi, d'o\`u  la d\'ecomposition sur $U$.
 
  On choisit un point g\'en\'eral $v$ de $V-U$. Il appartient  \`a une strate $S$ de $V$
 de dimension 
  strictement inf\'erieure \`a la dimension $n$ de $V$.
 La puret\'e locale en $v$ sert alors \`a \'etendre la d\'ecomposition au voisinage  du point $v$
dans une section normale $\NN_v$ dans $V$ \`a la strate $S$. Le point $v$ est dans une strate induite de dimension minimale sur $\NN_v$, ce qui nous ram\`ene \`a \'etudier le cas d'une strate de dimension minimale.
 On  suppose par construction que l'image r\'eciproque de la strate $S$ est un DCN dans $X$ qui est 
 relatif sur la strate $S$, on pourra \'etendre la d\'ecomposition sur tout un voisinage de $v$ dans $V$. 
 
 {\it Morphismes d'intersection}.  Avec notre hypoth\`ese sur $f$, il est naturel d'exprimer  le r\'esultat en termes de {\it morphismes d'intersection $I$} dont l'importance  apparait comme une des  r\'ev\'elations de la th\'eorie dans \cite{BBD}. Soit $X_{S_l}:= f^{-1}(S_l)$ l'image r\'eciproque d'une strate $S_l$ de dimension $l \leq n$,  o\`u $n$
  est la dimension de $V$, d'une stratification de Thom-Whitney $\SS$,  on pose $i_{X_{S_l}}: X_{S_l} \to X$, $i_{S_l}: S_l \to V$ et on consid\`ere le morphisme compos\'e:   
 \begin{equation}
 I_{S_l}:  R i_{X_{S_l}}^{!}\ilm \rightarrow \ilm
\rightarrow i_{X_{S_l}}^* \ilm 
\end{equation}
Sous l'hypoth\`ese de fibration sur les strates de $V$,  on en d\'eduit des  syst\`emes locaux 
images  
 \begin{equation}
 {\LL}^i_l = \im [{ R}^{-l+i} {f_l}_*
 ( R i_{X_{S_l}}^{!}   \ilm) \stackrel{I_{S_l}}{\rightarrow}
 R^{-l+i} {f_l}_* ( i_{X_{S_l}}^* \ilm)], 
  \end{equation}
sur les
diff\'erentes  strates  qui interviennent dans  la formule explicite de la d\'ecomposition de $ K := {R f}_{ *} \ilm $ en tant que  somme directe  d'extensions interm\'ediaires par  $i_{S_l}$:
$  \ph{i}{K} \simeq
\oplus_{S_l \in \SS }{i_{S_l}}_{! *} {\mathcal L}^i_l[l]$.  
Ces  ${\mathcal L}^i_l$  sont en  fait des VSH polaris\'ees de poids
$a+i-l$, car ${\LL}^i_l$  est l'image d'une variation de SHM
de poids $\omega \geq a+i-l$ dans une  variation de SHM
de poids $\omega \leq a+i-l$ et de plus le tout est calcul\'e avec des complexes logarithmiques 
en $X_{S_l}$, DCN relatifs sur $S_l$.
  
Retenons que le  th\'eor\`eme  de  d\'ecomposition sera  d\'emontr\'e par induction   sur la  dimension d\'ecroissante des  strates de $V$. 
La preuve n'utilise pas la th\'eorie des cycles \'evanescents et diff\`ere des preuves actuelles que nous citons en r\'ef\'erence  \cite{BBD, MI}.

\section{D\'ecomposition de la cohomologie perverse}
 Nous allons expliquer  le pas de la r\'ecurrence dans le cas d'une fibration
 par des DCN sur les strates ( \cite{EL},  d\'efinition ?). 
 La  situation est donc celle 
 d'un morphisme projectif de  vari\'et\'es alg\'ebriques complexes $f:X \to V$ o\`u $X$ est lisse
 et d'une  $VSH$  polaris\'ee d\'efinie  sur un  syst\`eme local
  $\tilde \LL$ sur  le compl\'ementaire de $Y$ un DCN dans $X$.
  Une stratification   $\SS$ de $V$ d\'efinit une famille de sous-espaces ferm\'es
 $V_i$, r\'eunion de strates de dimension $\leq i$ pour $i \in [0,n]$: $V_0 \subset  V_1 \subset \cdots \subset V_{n-1} \subset V = V_n$, 
et  $V_{n-1}$  contient le lieu singulier de $V$. Soit   $X_{V_i} = f^{-1}(V_i)$
et pour chaque strate $S_l$ de dimension $l$ soit  $X_{S_l} = f^{-1}(S_l)$,  par  construction, l'espace $X_{V_i} $ est un DCN dans $X$. Les propri\'et\'es de la stratification que nous utilisons sont:
\begin{itemize}
\item (T) Au-dessus de $S_l$  le  morphisme $f_l: X_{S_l} \to S_l$, induit par $f$,
 est un fibr\'e topologique localement trivial de fibre en un point $v$  un DCN dans l'image r\'eciproque lisse d'une setion normale g\'en\'erale $\NN_v$  \`a $S_l$ en $v$.
 \item (W) Le link  de la strate $S_l$ en un point de $S_l$ est un invariant topologique localement constant.
\end{itemize}
 Comme les stratifications sont de Thom-Whitney, le r\'esultat de J. Mather dans \cite{Ma}
nous garantit la propi\'et\'e W.

  En particulier, la  restriction de $f$ \`a 
$X - X_{V_{n-1}}$ est lisse sur  $ V-  V_{n-1} $. 
On peut toujours supposer  $Y = X_{V_{n-1}} \cup Y_h$ union de $X_{V_{n-1}} $ avec un DCN horizontal 
  $Y_h$ tel que la restriction  de  $f$ \`a $Y_h $ induise  au-dessus de de la grande strate 
 $U :=  V-  V_{n-1}$ un DCN relatif.  De m\^eme la fibre de $ X_{S_l}$ en un point $v$ est  un DCN dans l'image r\'eciproque d'une section normale g\'en\'erale $\NN_v$ \`a $S_l$ en $v$.
Dans ce  cas la stratification est dite adapt\'ee \`a $f$ et $\LL$.

  L'espace
$Y$  contient donc les singularit\'es de $\LL$ et de  $f$ de sorte que l'on puisse utiliser  des complexes  logarithmiques pour  la th\'eorie de Hodge  \`a coefficients dans $\LL$ sur $X$  \cite{EII}.
 Le diagramme suivant r\'esume les notations
   $$\begin{array}{ccccc}
 X - X_{V_l}  &
\stackrel{j_l}{\rightarrow}& X & \stackrel{i'_l}
{\hookleftarrow}&X_{V_l}\\
\downarrow&\quad&\downarrow f &\quad&\downarrow f_l  \\
 V- V_l  &\stackrel{k_l}{\rightarrow}&
   V & \stackrel{i_l}{\hookleftarrow}&
V_l\end{array}$$
 o\`u  $i'_l: X_{V_l} = f^{-1}(V_l) \to X, \,i_l: V_l  \to V $, (resp.$j_l: (X - X_{V_l}) \to X, \, k_l : (V - V_l) \to V$) 
 d\'enotent les plongements dans $X$ et $V$, (resp. les plongements des compl\'ements ouverts). 
 Soit $v$  un point dans $V_0$, on \'ecrit  $i_{X_v}: X_v = f^{-1}(v) \to X, 
 \, k_v : (V-\{v\}) \to V $ et $j_v: (X - X_v)\to X$. Enfin    $V_l^*:= V_l - V_{l-1}$ est la r\'eunion (\'eventuellement vide) des strates  de dimension $l$.
\\

On va montrer que la puret\'e  locale en $v$ permet  d'\'etendre la d\'ecomposition de $V- V_0$ \`a $V$ 
 \`a travers les points $v \in V_0$ de la r\'eunion des strates de  dimension z\'ero.
Cet argument s'applique en fait \`a l'\'etape  de r\'ecurrence  inductive   \`a travers une
strate de dimension quelconque, auquel cas le point $v$ est l'intersection avec une
section   normale  \`a  la strate $\NN_v$ dans $V$.  On \'ecrit aussi  pour chaque strate:
 $i_{S_l}: S_l   \to \overline S_l$ et de m\^eme avec abus  de notation  $i_{S_l}: S_l   \to V$. 
 
 \begin{prop}[D\'ecomposition de la cohomologie perverse] \label{P}
 Soit $i_0:V_0 \to V$ la r\'eunion des  strates de dimension $0$ et supposons 
 que la cohomologie perverse  de  $K = {R f}_{ *}  \ilm$ se  d\'ecompose sur l'ouvert  $ k_0: (V - V_0) \to V$ en une somme directe  d'extensions  interm\'ediaires   de VSH : ${\LL}^i_l$ (1.2) sur 
  toutes les strates $S_l$ de $V_l^*:= V_l-V_{l-1}$ pour tout $l > 0$  et en tout degr\'e $i$, soit:
\begin{equation}
 {^p\mathcal { H}}^i
( K)_{\vert V-V_0} \simeq
 \oplus^{S_l \subset V_l^*}_{0 < l \leq n } \;\;  k_{0}^* {i_{S_l}}_{! *}{\LL}^i_l[l], \quad K_{\vert V-V_0}  =   \oplus_{i \in \Z}{^p\mathcal { H}}^i
( K)_{\vert V-V_0} [-i] 
\end{equation}

 Alors la suite exacte longue de cohomologie perverse 
 \begin{equation*}
{^p\mathcal { H}}^i ((i_0)_* R(i_0)^! K)
\stackrel{^p\!\alpha_i}{\rightarrow}{^p\!\mathcal { H}}^i ( K)
     \stackrel{^p\!\rho_i}{\rightarrow}{^p\!\mathcal { H}}^i ( R{k_0}_*
      K_{\vert V - V_0})
    \stackrel{^p\!\delta_i}{\rightarrow}
\end{equation*}
donne  lieu \`a une suite exacte courte: 
$0 \to  Im\,  {^p\!\alpha_i}  \to {^p\!\mathcal { H}}^i ( K)  \xrightarrow{^p\!\rho_i}  Im \,{^p\!\rho_i} \to 0$ 
de faisceaux pervers, scind\'ee sur $V$, qui  se d\'ecompose en termes de  ${\LL}^i_l$ et  ${\LL}^i_0$ (1.2):

\smallskip
$ Im \,{^p\!\rho_i} =  (k_0)_{! *} k_0^*\, \ph{i}{ K} \simeq
\oplus^{S_l \in V_l^*}_{0 < l \leq n } {i_{S_l}}_{! *} {\LL}^i_l[l]$, \,et $ \,  Im\,  {^p\!\alpha_i} =  ker \, {^p\!\rho_i} \simeq (i_0)_* {\LL}^i_0. $
 \end{prop}
 
\noindent La preuve \'etant locale aux points de $V_0$, on peut supposer $ V$ projective et $V_0$
   r\'eduit \`a un point $v$  et l'on \'ecrit  $i_v, k_v$ pour les immersions au lieu de  $i_0, k_0$. La suite exacte dans la proposition est associ\'ee au triangle sur $V$: $
 i_{v *} R{i}_v^{!} (K)
   \stackrel{\alpha}{\rightarrow}   K
   \stackrel{\rho}{\rightarrow} {R k}_{v *} K_{\vert V-\{v\}}
    \stackrel{[1]}{\rightarrow }$,  et 
    de plus on a:
 ${^p\mathcal { H}}^i (i_{v *} R{i}_v^{!} K) \simeq 
i_{v*} H^i ( R{i}_v^{!} K) $.
Afin de  calculer  successivement: $Im \, {^p\!\rho_i}$,  $ Im\,  {^p\!\alpha_i}$ et prouver le scindage dans $\ph{i}{ K}$, 
il nous est utile d'abord  de signaler le calcul  suivant de cohomologie perverse de $ R k_{v *} K_{\vert V-\{v\}}$ \`a partir de la d\'ecomposition (2.1): 

\begin{lemma} 
1) Dans le cas d'un seul syst\`eme local $\LL'$ sur une strate $S_l$ contenant $v$ dans son adh\' erence, on note  $i_{S_l}: S_l   \to \overline S_l,\, k_v:(\overline S_l-v)
 \to \overline S_l$,  alors  la suite longue de cohomologie perverse d\'efinie par le triangle:
$$ i_v^! (i_{S_l})_{ !*}{\LL}'[l] \to (i_{S_l})_{ !*}{\LL}'[l] \to
Rk_{v *}k_v^*(i_{S_l})_{ !*}{\LL}'[l])\xrightarrow{[1]}$$
  se d\'ecompose en une suite exacte 
\begin{equation*}
0 \to (i_{S_l})_{ !*}{\LL}'[l] \to {^p\!\mathcal { H}}^0 (Rk_{v
*}k_v^*(i_{S_l})_{ !*}{\LL}'[l]) \to i_{v,*} H^1_v ((i_{S_l})_{ !*}{\LL}'[l])
  \to 0
\end{equation*}
de plus, on a des isomorphismes:
 \begin{equation*}
 \begin{split} H^0({^p\mathcal H}^0 (Rk_{v *} k_v^*(i_{S_l})_{ !*}{
  \LL}'[l])) & \simeq  R^0 k_{v *}k_v^*(i_{S_l})_{ !*}{\LL}'[l]\simeq i_{v,*} H^1_v
   ((i_{S_l})_{ !*}{\LL}'[l]) \\
    R^i k_{v *}k_v^*(i_{S_l})_{ !*}{\LL}'[l] &\simeq i_{v,*} H^{i+1}_v
   ((i_{S_l})_{ !*}{\LL}'[l]) \hbox {\, \rm pour \, } i > 0 \\
     {^p\mathcal H}^i (R k_{v *} k_v^*(i_{S_l})_{ !*}{\mathcal L}'[l])
& = 0 \hbox {\rm \, pour \,} i < 0,\\
 {^p\mathcal H}^i (Rk_{v *} k_v^*(i_{S_l})_{ !*}
 {\mathcal L}'[l]) & \simeq H^{i+1}_v ((i_{S_l})_{ !*}{\LL}'[l]) \hbox {\, \rm pour \, } i > 0.
 \end{split}
\end{equation*}
 2) Dans le cas d'une somme directe de  syst\`emes locaux  $K' := \oplus_j (i_{S_l})_{ !*}{\mathcal L}_l^j
[l-j]$  sur la m\^eme strate $S_l$, on a: \\
 i)  Une suite exacte  courte pour tout $i$
\begin{equation*}
 0 \to
(i_{S_l})_{ !*}({\mathcal L}_l^{i}[l]) \to {^p\mathcal H}^{i}(R
k_{v *}k_v^* K') \xrightarrow{h} \oplus_{j \leq i} R^i k_{v *}
(k_v^*(i_{S_l})_{ !*}{\mathcal L}_l^j [l-j]) \to 0
 \end{equation*}
o\`u le dernier terme est un faisceau en degr\'e z\'ero de support $v$.\\
\n ii) $H^0( i^*_v {^p\mathcal { H}}^{i}(R k_{v *}k_v^* K'))
\simeq \oplus_{j \leq i} R^i k_{v *} ( k_v^*(i_{S_l})_{ !*}{\mathcal
L}_l^j [l-j]) \simeq R^i k_{v *}(k_v^* ({^p\!\tau_{\leq i}}K')) $\\
\n iii)  En particulier un morphisme $\varphi$ \`a valeur  dans ${^p\!\mathcal
H}^{i}(R k_{v *}k_v^* K') $ se factorise par \\ $(i_{S_l})_{ !*}
({\mathcal L}_l^{i}[l])$ si $h \circ \varphi = 0$.
  \end{lemma}
En 1,  on utilise:  ${^p\!\mathcal { H}}^i ((i_{S_l})_{ !*}{\LL}'[l]) = 0 $
pour $i \not =  0$, $H^0(i_v^*(i_{S_l})_{ !*}{\LL}'[l]) = 0$ et
$H^i_v((i_{S_l})_{ !*}{\LL}'[l]) = 0$ pour $ i \leq 0$.

En  2 i) on applique  l'assertion 1  \`a
${\LL}':= {\mathcal L}_l^j$ sur $S_l$ pour  d\'eduire:\\
${^p\mathcal { H}}^i (Rk_{v *}k_v^* (i_{S_l})_{ !*}{\mathcal
L}_l^j[l-j]) = {^p\mathcal { H}}^{i-j} (Rk_{v *}k_v^* (i_{S_l})_{ !*}{\mathcal
L}_l^j[l]) = 0 $ pour $i < j$,  \\
${^p\mathcal { H}}^i (R k_{v *}k_v^* (i_{S_l})_{ !*}{\mathcal
L}_l^j[l-j]) = R^i k_{v *}
 (k_v^* (i_{S_l})_{ !*}{\LL}_l^j[l-j]) $ pour $i > j$,\\
   et la suite  exacte pour $ i = j$: \\
$ 0 \to  (i_{S_l})_{ !*}{\mathcal L}_l^j[l] \to {^p\mathcal {
H}}^{j} (R k_{v *}k_v^* (i_{S_l})_{ !*}{\LL}_l^j [l-j]) \to R^{j}
k_{v *}
 (k_v^* (i_{S_l})_{ !*}{\LL}_l^j [l-j]) \to 0.$\\
 L'\'enonc\'e 2 ii) se d\'eduit du cas $i < j$, et 2 iii) s'applique pour toute suite exacte dans une cat\'egorie ab\'elienne.

 \subsection{ Preuve de la proposition \ref{P}}  1) {\it Calcul de l'image de  $\, {^p\!\rho_i}$}.
 Par d\'efinition  du foncteur extension interm\'ediaire \cite{BBD}, 1.4.22, p.54, on a:  
 $$  (k_{v})_{!*}\ph {i}{k_v^*K}  = Im\, (\ph {i}{R
(k_{v})_! k_v^*K}  \xrightarrow{^p\!\gamma_i} \ph {i}{R (k_{v})_* k_v^*K})$$
qui est \'egal \`a $\oplus^{S_l \subset  V_l^*}_{0 < l \leq n }(i_{S_l})_{! *} {\LL}^i_l[l] $, $V_l^* = V_l - V_{l-1}$, dim.$S_l = l$,  d'apr\`es l'hypoth\`ese (2.1).
 Le morphisme ${^p\!\gamma_i}$ se  factorise en:
 \begin{equation*}
 \ph {i}{R (k_{v})_! k_v^*K} \xrightarrow{^p\!\beta_i} \ph {i}{K}
\xrightarrow{^p\!\rho_i} \ph {i}{R (k_{v})_* k_v^*K},  \quad {^p\!\gamma_i} = {^p\!\rho_i}\circ {^p\!\beta_i} 
\end{equation*}
 on en d\'eduit:
\begin{equation*}
 Im\, {^p\!\rho_i} \supset  \oplus^{S_l \in V_l^*}_{0 < l \leq n } (i_{S_l})_{! *} {\mathcal L}^i_l[l] = Im\, {^p\!\gamma_i}.
\end{equation*}
Pour obtenir  l'\'egalit\'e  $Im\, {^p\!\rho_i}  = Im\, {^p\!\gamma_i}$,
il suffit de prouver, d'apr\`es le lemme 2.2  (2~iii),  que le morphisme  $ {^p\!\rho_i^0}$,  
induit sur la  cohomologie en degr\'e z\'ero par  $ {^p\!\rho_i}$  s'annule:
\begin{equation*}
  H^0(i^*_v\,\ph {i}{ K})  \stackrel{ {^p\!\rho_i^0}}{ \longrightarrow }
    H^0(i_v^*
\,\ph {i}{R k_{v *} k_v^*K}) \simeq \oplus^{S_l \in V_l^*}_{0 < l \leq n, j \leq i } 
R^i k_{v *} (k_v^* {i_{S_l}}_{! *}{\LL}_l^j [l-j]).
\end{equation*}
{\it Interpr\`etation du morphisme $ {^p\!\rho_i^0}$.}
 On consid\`ere de petits voisinages $B_v $  de $v$ et $B_{X_v}$ de $X_v$ et on interpr\`ete
le morphisme $ {^p\!\rho_i^0}$ comme un morphisme de 
$ \H^0(B_v, \,\ph {i}{ K}) $ dans \\ $ H^0(i^*_v\,\ph {i}{R k_{v *} k_v^*K})\simeq \H^{0}
(B_v - v, \ph {i}{K}) \simeq Gr^{\pt}_i\H^{i} (B_{X_v} - X_v, \ilm
)$. \\
En effet:   $  \H^i(B_v, \,{^p\!\tau}_ {\leq i}  K) \simeq  \H^i(B_v, \ph {i}{ K}[-i]) $
car  $ \H^i(B_v, \,{^p\!\tau}_ {< i} K) = H^i(i_v^* ({^p\!\tau}_ {< i} K)) = 0 $ et 
$ \H^{i+1}(B_v,\,{^p\!\tau}_ {< i}  K) = H^{i+1}(i_v^* ({^p\!\tau}_ {< i} K)) = 0 $, et par cons\'equent
 $ {^p\!\rho_i^0}$
  se factorise comme suit
 \[
 \H^i(B_v, \,{^p\!\tau}_ {\leq i}  K) \to  \H^i(B_v,  K) 
=  \H^{i} (X_v,  \ilm) \to \H^{i} (B_{X_v} - X_v,  \ilm )
  \]
 o\`u l'espace $ \H^{i} (X_v,  \ilm)$ est de poids  $\omega \leq a+i$ alors qu'\`a droite le poids est 
  $\omega > a+i$ d'apr\`es le r\'esultat sur   la semi-puret\'e en $v$, donc $ {^p\!\rho_i^0} = 0$.

\smallskip
 \n 2)  {\em L'isomorphisme $Im \,{^p\!\alpha_i}\simeq i_{v *}{\LL}^i_0$}.  
     Le morphisme de connexion 
     \begin{equation*}
{^p\!\mathcal { H}}^{i-1} ( {{R k}_{v}}_*   K_{\vert V-\{v\}})
 \stackrel{^p\!\delta_{i-1}}{\rightarrow} \H^i_v(V, K)
\end{equation*}
s'annule sur l'image de $ \,{^p\rho_{i-1}}$ et par cons\'equent,  d'apr\`es le lemme  pr\'ec\'edent 2.2,  (2 i), 
son image est \'egale  \`a celle de  l'espace vectoriel $R^{i-1} k_{v *}k_v^* ({^p\!\tau_{\leq i-1}}K)$   \begin{equation*}
R^{i-1} k_{v *}k_v^* ({^p\!\tau_{\leq i-1}}K) 
 \stackrel{ ^p\!\delta_{i-1}}{\rightarrow}{^p\!\HH}^{i}(i_v^!K) =  \H^i_v(V, K) \simeq \H^i_{X_v} (X,\ilm)
\end{equation*}
 et l'on a $Im \,{^p\!\alpha_i}\simeq \H^i_v(V, K)/ Im { ^p\!\delta_{i-1}}$. Par  ailleurs, rappelons que ${\LL}^i_0$
est l' image du morphisme d'intersection    $I_v^i$  en degr\'e $i$ dans le diagramme suivant
\begin{equation*}
H^{i-1} ({i^*_v}{ R k_{v}}_* K_{\vert V-\{v\}} )
\stackrel{\delta_{i-1}}{\rightarrow}
 \H^i_v (V, K) \stackrel{I_v^i}{\rightarrow} H^i({i^*_v}  K)
  \stackrel{\rho_i}{\rightarrow}
   H^i ({i^*_v} R k_{v*} K_{\vert V-\{v\}} )
\end{equation*}
et l'on a $Im\,  I_v^i \simeq \H^i_v(V, K)/ Im \, { \delta_{i-1}}$. Il suffit donc de prouver: $Im\,  {^p\!\delta_{i-1}} = Im\, {\delta_{i-1}}$.
 Vu que la d\'ecomposition s'applique  par r\'ecurrence sur $V-\{v\}$,
on trouve:\\
 $ {^p\!\delta_{i-1}} ( R^{i-1} k_{v
*} k_v^* ({^p\!\tau_{\leq i-1}}K)) \simeq
\delta_{i-1}({^p\!\tau_{\leq i-1}}\H^{i-1}(B_{X_v}-X_v, \ilm))$.
Alors que l'on veut l'\'egalit\'e avec  toute l'image $ \delta_{i-1}(\H^{i-1}(B_{X_v}-X_v, \ilm)) $,
ce qui d\'ecoule de la semi-puret\'e:
 en effet, le  quotient 
$\H^{i-1}(B_{X_v}-X_v, \ilm))/ {^p\!\tau_{\leq i-1}}$ est de poids
$\omega < a+i$  et le poids
de  $\H^i_{X_v} (X, \ilm)$ est $\omega \geq a+i$
d'o\`u les images   de $
{^p\!\delta_{i-1}} $ et $ \delta_{i-1}$ sont \'egaux dans $\H^i_{X_v} (X, \ilm)$
de poids $\omega \geq a+i$. Vu que  $Im \,{^p\!\alpha_i} = ker \, \,{^p\rho_i}$,  on obtient 
une suite exacte
\[
0 \to i_{v *}{\LL}^i_0 \to {^p\!\HH}^{i}(K) \to
\oplus^{S_l \subset  V_l^*}_{0 < l \leq n } {i_{S_l}}_{! *} {\LL}^j_l[l] \to 0
 \]
3) {\it  Scindage de $ {^p\!\HH}^{i}(K)$}. Consid\'erons  
la suite exacte: $ \ph {i}{R (k_{v})_! k_v^*K}
\stackrel{ ^p\!\beta_i}{\rightarrow}   \ph {i}{K} \stackrel{\theta_i}{ \rightarrow} H^i({i^*_v} K)
$. Par  d\'efinition de $ {\LL}^i_0 $, le morphisme 
$\theta_i$ induit un isomorphisme  sur $i_{v *}{\LL}^i_0 = Im \,{^p\!\alpha_i}
 = ker \,  {^p\!\rho_i}$, alors que $\theta_i \circ  {^p\!\beta_i} = 0$,
 et par cons\'equent
 $Im \, {^p\!\beta_i} \cap {\LL}^i_0 = 0$. On en d\'eduit que
 ${^p\!\rho_i}$  induit  un  isomorphisme:
  $Im \,{^p\!\beta_i}\xrightarrow{^p\!\rho_i} \, Im\,
{^p\!\rho_i} $ et finalement: 
 $  \ph {i}{K} = Im \, {^p\!\beta_i} \oplus i_{v *} {\LL}^i_0 $ 
\begin{remark} \label{r} On retient les relations utiles pour la suite
\begin{equation*}
\begin{split}
ker I_v^i = ker \ {^p\alpha_i} \simeq Im \ {^p\delta_{i-1}} &\simeq
\oplus^{S_l \subset  V_l^*}_{0<l \leq n, 0 \leq i-1-j } R^{i-1-j} k_{v *} ({i_{S_l}}_{!
*} {\LL}_l^j [l])\\
&\simeq \oplus^{S_l \subset  V_l^*}_{0<l \leq n, 0 < i-j  } H^{i-j}_v ({i_{S_l}}_{! *}
{\LL}^j_l[l]).
\end{split}
\end{equation*}
  \end{remark}
 \subsection{Lefschetz difficile}  Pour terminer, il faut d\'emontrer  l'isomorphisme de Lefschetz pour le cup-produit 
it\'er\'e  avec  la classe  $ \eta$ d'une section hyperplane sur la  cohomologie perverse.
Par r\'ecurrence,  il reste \`a v\'erifier le cas   d'un  point $v$ de la strate de  dimension z\'ero.
Nous le v\'erifions  sur les termes de la  formule   explicite  de  d\'ecomposition.
Donc  on  suppose   Lefschetz  difficile  pour
 $ (k_v)_{!*} k_v^*\ph{i}{K}$, et  on pose
 $ {\LL}^i_v = \oplus_{ v\in V_0}Im (\H^{i}_{X_v}
(X, \ilm) \stackrel{I_v^i}{\rightarrow} \H^i(X_v, \ilm))$ pour l'image du morphisme d'intersection $I_v^i$. Il reste \`a  prouver
\begin{prop} i) 
 La $SH$ ${\LL}^{-i}_v$ est  duale de  Poincar\'e de
${\LL}^{i}_v$  pour tout $i \in \Z$.\\
ii) Le cup-produit it\'er\'e avec $\eta$ 
induit des 
isomorphismes  $\eta^i:{\LL}^{-i}_v \to {\LL}^{i}_v$ pour $i  \geq 0$.
En particulier,  ${\LL}^{i}_v$ est une SH   polarisable.
 \end{prop}
\n  L'assertion i) est claire sur la d\'efinition auto-duale du morphisme d'intersection. \\
ii) {\it Interpr\'etation  dans le cas classique} o\`u $ \tilde \LL$ est un syst\`eme local  en degr\'e $0$ sur $X$ lisse et  $X_v$ aussi est lisse. Le r\'esultat   utilise  la d\'ecomposition de Lefschetz sur $X_v$ comme suit
$$
\begin{array}{ccc}
 \H^{-(i-1)} (X_v,\tilde \LL_{\vert X_v} [m-1]) &  \overset { \eta^{i-1}} {\overset { \simeq}{\longrightarrow}}  &  \H^{(i-1)} (X_v,\tilde \LL_{\vert X_v} [m-1])  \\
 &&\\
 \uparrow \eta &  \searrow \eta^i & \downarrow \eta  \\
 &&\\
   \H^{-(i+1)} (X_v,\tilde \LL_{\vert X_v} [m-1])&   \overset  { \eta^{i+1}}{\overset { \simeq}{\longrightarrow}}    & 
   \H^{(i+1)} (X_v,\tilde \LL_{\vert X_v} [m-1])
\end{array}
$$
Pour interpr\'eter ce diagramme, on note
que pour $\LL:= \tilde \LL [m] $, l'image de l'injection $\eta$ \`a gauche est \'egale   
\`a l'image   ${\LL}^{i}_v := Im \, ( \H^{-i}_{X_v}(X, \LL) \xrightarrow{I^{i}_v} 
  \H^{-i} (X_v, \LL_{\vert X_v}))$ d'apr\`es  l'isomorphisme de Thom-Gysin:
  $\H^{-(i+1)} (X_v,\tilde \LL_{\vert X_v} [m-1]) =  \H^{-i-2}(X_v, \ilm) \simeq  \H^{-i}_{X_v}(X, \ilm)$, alors que  l'image de la surjection $\eta$ \`a droite est \'egale \`a  l'image ${\LL}^{i}_v$ de $   \H^i_{X_v}(X, \LL)  \stackrel{I_v^i}{\rightarrow} \H^i(X_v, \ilm))$.
 
 Donc,  l'isomorphisme $ \eta^i$, induit  sur l'image de l'injection $\eta$ \`a gauche, un 
isomorphisme 
 sur l'image de  $\eta$ \`a droite:
 $${\LL}^{-i}_v \simeq \eta( \H^{-(i+1)} (X_v,\tilde \LL_{\vert X_v} [m-1])) \overset{  \eta^i}{ \overset{\simeq}{\longrightarrow}}  \H^{i+1} (X_v,\tilde \LL_{\vert X_v} [m-1]) \simeq {\LL}^i_v$$ qui se  d\'eduit  de  l'isomorphisme classique $ \eta^{i+1}$ en bas
 dans le diagramme. 

  La partie primitive de $\H^{-i} (X_v, \LL_{\vert X_v})$ n'est donc pas concern\'ee, ce qui pourrait expliquer le r\'esultat surprenant qui suit et qui \'evite le cas crucial dans le cas classique d'une r\'ecurrence sur une section hyperplane.\\
  {\it Preuve de ii).}  L'\'enonc\'e garde un sens m\^eme si $X_v$ n'est pas la fibre d'un morphisme, mais on va  utiliser  la remarque (2.3) pr\'ec\'edente dans la preuve
 et donc le fait que $X_v$ soit la  fibre d'un morphisme $f$. 
 On proc\`ede  par  r\'eduction \`a une section  hyperplane relative  lisse $H$ dans $X$,
  d'intersection transversale aux strates du DCN: $X_v \cup Y$. Soit   $i_H: H \to X$ l'immersion. Le cup-produit 
   avec la  classe  $\eta$ de  $H$  d\'efinit un  morphisme  \'egal  au  compos\'e des
  morphismes $ \ilm  \stackrel{ \rho}{\rightarrow} i_{H*}i^{*}_H \ilm \stackrel{G}{\rightarrow}  \ilm
[ 2]$. Par application des foncteurs $R i^{!}_{X_v}$ et
$i^{*}_{X_v}$ on obtient des morphismes  commutant avec  les  morphismes d'intersection $I_v^*$ sur $X$ et $I_v^*(H)$ sur $H$  dans le diagramme
$$\begin{array}{ccccc} \H^i_{X_v}(X, \ilm)
&\stackrel{ \rho ^!}{\rightarrow}&\H^i_{X_v \cap H} (H, \ilm)
 &\stackrel{ G^!}{\rightarrow}& \H^{i+2}_{X_v}(X, \ilm)\\
I_v^i{\downarrow}\quad  &   &I_v^i (H){\downarrow}\quad  && I_v{i+2}{\downarrow}\qquad  \\
\H^i(X_v, \ilm)
 &\stackrel{ \rho^*}{\rightarrow}&\H^i ({X_v \cap H} ,\ilm)
 &\stackrel{G^*}{\rightarrow}&  \H^{i+2}(X_v, \ilm)
 \end{array}$$
On pose:
${\LL}^i_v = Im \, I_v^i,\, {\LL}(H)^i_v = Im\,  I_v^i(H),\,
  {\LL}^{i+2}_v = Im\,  I_v{i+2} $ pour les  images des morphismes  verticaux
  d\'efinis par $I_v^*$. On s'int\'eresse \`a  l'image de la premi\`ere ligne dans la seconde ligne repr\'esent\'ee  par   la ligne interm\'ediaire qui aurait d\^u figurer dans le diagramme:
$$
 {\LL}^i_v \stackrel{ \rho'}{\rightarrow}
 {\mathcal L}( H)^i_v \stackrel{ G'}{\rightarrow}
{\LL}^{i+2}_v $$
o\`u $\rho'$ ( resp. $G'$ ) est induit par le morphisme $ \rho^*$
    (resp. $ G^* $ ). 
  \begin{lemma} Le morphisme induit $\rho':{\LL}^i_v \to {\LL}(H)^i_v$
 est un isomorphisme pour $i < 0$
et par dualit\'e $G': {\LL}(H)^i_v \to
  {\LL}^{i+2}_v $ est un isomorphisme pour $i > 0 $.
 \end{lemma}
    \n {\it  Preuve}. D'apr\`es le th\'eor\`eme d'annulation d'Artin-Lefschetz
    sur l'espace affine $X_v - X_v\cap H$, on a:
\begin{lemma} i) Les morphismes $  \rho ^!:
\H^i_{X_v}(X, \ilm) \rightarrow \H^i_{X_v\cap H} (H,
\ilm) $ \\ sont des isomorphismes pour $i <0$  et injectifs pour $i = 0$.\\
ii) Les morphismes $\rho ^*: \H^i ({X_v}, \ilm) \rightarrow \H^i
({X_v \cap H} , \ilm) $\\ sont des isomorphismes pour $i < -2$ et 
injectifs pour $i= -2$.
 \end{lemma}
\n  On d\'emontre  i). Le complexe
 $P_{X_v} = R i^{!}_{X_v}  j_{!*}{\LL}[1]$ est
 un faisceau pervers sur le diviseur de Cartier  $X_v$ \cite{BBD} p.106.
Sur l'espace  affine, on a:
 $\H^j_c( X_v - X_v \cap H, P_{X_v}) = 0 $ pour $j < 0$, d'o\`u on d\'eduit  l'isomorphisme   pour  $i < 0 $\\
 \centerline {$ \H^i_{X_v}(X, \ilm) \simeq  \H^{i-1}(X_v, P_{X_v})
 \stackrel{ \rho ^!}{\rightarrow}
\H^{i-1} ({X_v \cap H}, P_{X_v}) \simeq \H^{i}_{X_v \cap H} (H, \ilm)$. } 

\vskip.2in

\n {\it Suite de la preuve.}  
On pose  $K = R f_* \ilm$ et  $K(H) = R(f_{\vert H})_*
 (\ilm_{\vert H})$  et, vue la remarque (\ref{r}),  on introduit les morphismes $\delta_{i-1}$ sur $X$, et $\delta_{i-1}( H) $
sur $H$, dans  le diagramme suivant:
$$\begin{array}{ccc}
 H^{i-1} (i^*_v R(k_v)_* K_{\vert V-\{v\}}) &
\stackrel{ \rho_v}{\rightarrow} & H^{i-1} ({i^*_v} R (k_v)_ * {K(H)}_{\vert V-\{v\}})\\
 \delta_{i-1} \downarrow && \delta_{i-1}( H) \downarrow \quad \\
\H^i_{X_v}(X, \ilm)\simeq \H^{i}_v (V,K)&\stackrel{ \rho^!}
{\rightarrow}&\H^i_{X_v \cap H} (H,\ilm)\simeq \H^i_v (V,K(H))\\
  I_v^i \downarrow &&
I_v^i(H) \downarrow \quad  \\
 \H^i(X_v, \ilm) \simeq \H^{i} (i^*_v K)
  &\stackrel{
\rho^*}{\rightarrow}&\H^i ({X_v \cap H} ,\ilm)
  \simeq \H^{i}(i_v^* K(H))
 \end{array}$$
Pour \'etablir l'isomorphisme 
  $\rho':{\LL}^i_v \simeq {\LL}(H)^i_v$,  \'etant donn\'e que $ \rho^!$ est un isomorphisme pour $i<0$, il suffit de prouver 
 que  le morphisme  $ \rho^!$ induit aussi un isomorphisme   $ker \, I_v^i \simeq  ker \, I_v^i(H)$ ou aussi
 $ Im \, \delta_{i-1} \simeq  Im \, \delta_{i-1}(H)$;  or on a  d'apr\`es la remarque pr\'ec\'edente (\ref{r}):
\begin{enumerate}
   \item $  Im \, \delta_{i-1}
= Im \,(\oplus^{S_l \subset  V_l^*}_{0 < l \leq n, j < i }
 H^{i-1-j} ( {i^*_v}
 {R k_v}_* k_v^*{i_{S_l}}_{! *}\LL^j_l [l]))$
   \item $  Im \, \delta_{i-1}(H) =
   Im (\oplus^{S_l \subset  V_l^*}_{0 < l \leq n, j < i }
 H^{i-1-j}(i^*_v R (k_v)_* k_v^*{i_{S_l}}_{! *}
 (\LL_{\vert H})_l^j [l])) $
 \end{enumerate}
 (noter que:  la restriction $\LL_{\vert H}[-1]$ d\'ecal\'ee de
$[-1]$ est un faisceau pervers sur $H$, et  d'apr\`es la formule (1.2): $(\LL_{\vert H}[-1])^{j+1}_l = (\LL_{\vert H})^j_l$ et que  $ I_v^i(H, \ilm)) = I_v^{i+1}(H, \ilm[-1])$), ce qui compense les indices).
Pour comparer  les faisceaux $\LL_l^j$ et $(\LL_{\vert H})_l^j$, on  consid\`ere un point  $u_l$ d'une strate  $S_l$
de  dimension  $l$,
  une section $N_{u_l}$   normale \`a $S_l$ au  point
$u_l$,  et les  sous-espaces $ X_{N_{u_l}}:= f^{-1}(N_{u_l})$  qui est lisse et $X_{u_l} := f^{-1}(u_l)$
qui est un  DCN dans $ X_{N_{u_l}}$  par construction, alors:
 \[
 (R^{-l+j} f_*(i_{X_l})^{! }{\ilm})_{u_l}\simeq
\H^{-l+j}_{X_{u_l}}(X_{N_{u_l}},{\ilm})\simeq
\H^{j}_{X_{u_l}}(X_{N_{u_l}},{\ilm}[-l])
 \]
 o\`u la restriction de ${\ilm}[-l]$  \`a $X_{N_{u_l}}$
est un faisceau pervers car $X_{N_{u_l}}$ est de  codimension  $l$ par transversalit\'e. Comme $i < 0$ et $ j
< i $ on a $ j < -1 $ et le th\'eor\`eme d'annulation d'Artin sur l'espace affine
 $X_{u_l}-(H\cap X_{u_l})$ s'applique;  
on en d\'eduit l'isomorphisme:  
$$R^{-l+j} f_*({i_{X_l}}^{! }{\ilm}) \simeq  R^{-l+j} 
f_*(({i_{(H \cap X_l})}^{! }{\ilm}_{\vert H})$$
des deux termes pris avant  restriction \`a $H$   \`a gauche  et apr\`es restriction \`a $H$   
\`a droite.
 En particulier le raisonnement s'applique aussi sur une strate g\'en\'erique $S_n$. En effet  $X_n := f^{-1}(S_n)$ est ouvert, et  en tout point 
 $u$,  l'espace $X_{u} := f^{-1}(u)$ est lisse, on a:
 \[
 (R^{-n+j} f_*(i_{X_n})^{! }{\ilm})_u = (R^{-n+j} f_*(i_{X_n})^* {\ilm})_u \simeq
\H^{-n+j} (X_u, {\ilm})
 \]
o\`u  $ \H^{-n+j} (X_u, {\ilm})  =  \H^{m-n+j} (X_u, {\tilde \LL})$ est isomorphe \`a \\ $\H^{-n+j} (X_u, {\tilde \LL}) = \H^{m-n+j} (X_u \cap H, {\ilm})$ par le th\'eor\`eme sur la section hyperplane de Lefschetz
 car   $X_u$ est de dimension  $m-n$ et $j < -1$.
D'o\`u  $\rho^!$  qui est un  isomorphisme, induit aussi des isomorphismes de $Im \, \delta_{i-1} = ker \, I_v^i $ sur $ Im \, \delta_{i-1}(H) = ker \, I_v^i(H)$. Les colonnes du diagramme  \'etant exactes,   il induit les  isomorphismes $\rho'$ pour $i > 0$  (et par dualit\'e on a des isomorphismes $G'$) dans le diagramme suivant:
 $$\eta^i: {\LL}^{-i}_v    \stackrel{ \rho'}{\simeq} {\LL}(H)^{-i}_v \simeq( {\LL}_{\vert H })^{-i+1}_v
\stackrel{ \eta^{i-1}}{\simeq} ({\LL}_{\vert H })_v^{i-1}
\simeq {\LL}(H)^{i-2}_v  \stackrel{ G'}{\simeq}
  {\LL}^i_v $$
o\`u $ \eta^{i-1}$ est un isomorphisme par r\'ecurrence; ce qui termine la  preuve du lemme, et celle de la   proposition.

 \begin{corollary} i) Le cup-product avec la classe d'une section   hyperplane
  d\'efinit  des isomorphismes $\eta^i:{^p\HH}^{-i}(K) \to
{^p\HH}^i(K) $ pour $i
\geq 0$.\\
\n ii) Le th\'eor\`eme de d\'ecomposition est vrai pour un morphisme fibr\'e par des DCN sur les strates, donc pour tout morphisme alg\'ebrique propre.
 \end{corollary}
 \n ii)  En effet,  le th\'eor\`eme s'applique sur la grande strate d'apr\`es le cas classique.  Si on suppose par r\'ecurrence le th\'eor\`eme vrai sur $V- V_l$,  
   il se prolonge \`a la r\'eunion des strates $S_l$ de dimension $l$, car en tout point $v$ de $S_l$, les r\'esultats pr\'ec\'edents s'appliquent sur la normale en $v$ \`a $S_l$ ce qui entra\^ine le r\'esultat au voisinage de $v$ dans $S_l$ et donc sur tout $S_l$. On a vu que pour tout morphisme $f$, on peut se r\'eduire au cas fibr\'e par des DCN relatifs sur les strates.
 
  \begin{remark}[La d\'ecomposition naturelle]   
 Deligne d\'eduit du  th\'eor\`eme \\ de d\'ecomposition, pour  $\eta$ fix\'e, un morphisme canonique  $ \oplus {^p{\HH}}^i(K)[-i]  \stackrel{\gamma}{\rightarrow} K$  \cite{can} qui pr\'esente
l'avantage d'\^etre compatible avec la SHM sur la cohomologie d'un ouvert de la forme $f^{-1}(U)$,
ce qui devrait permettre de construire une th\'eorie de Hodge  en bas sur $V$ \`a partir des 
 VSH  sur les images $\LL^i_l$ des morphismes d'intersection et compatible avec la SHM induite sur les gradu\'es pervers utilis\'es dans cette note. Ce r\'esultat  correspond au fait 
 bien  connu  que la suite spectrale de Leray perverse est compatible avec les SHM. 
  \end{remark}  
\bibliographystyle{amsalpha}

\bibliographystyle{amsalpha}

\end{document}